\documentclass[12pt,leqno]{article}
\title{Some congruences on prime factors of class number of  finite algebraic extensions $K/\Q$,  version 4.0}
\author{Roland Qu\^eme}

\usepackage{amsmath,amsthm,amstext,amsbsy}
\newtheorem{thm}{Theorem}[section]
\newtheorem{cor}[thm]{Corollary}

\font\mathbb=msbm10
\newcommand{\N}{\mbox{\mathbb N}}

\newcommand{\Q}{\mbox{\mathbb Q}}

\newcommand{\modu}{\ \mbox{mod}\ }
\newcommand{\be}{\begin{equation}}
\newcommand{\ee}{\end{equation}}
\date{2002 april  20}
\begin{document}
%%% ====================================================================
\maketitle
\tableofcontents
\clearpage
\abstract
Roland Qu\^eme

13 avenue du ch\^ateau d'eau

31490 Brax

France

2003 april 20

tel: 0561067020

mailto: roland.queme@free.fr

home page: http://roland.queme.free.fr/

****************************


Primary  11R29:  Class Number

Secondary 12F10: Galois Theory

****************************

\begin{itemize}
\item 
This paper is a contribution to the description of some congruences on the odd prime factors of the class number of the number fields.
\item
We say that a finite Galois extension $L/K$ is  Galois solvable  if the Galois group
$Gal(L/K)$ is solvable. 
An  example  of the results obtained  is:

{\it 
Let $L/\Q$ be a   finite Galois solvable   extension   with $[L:\Q]= N$, where $N>1$ is odd.
Let $h(L)$ be the class number of $L$. Suppose that $h(L)>1$.
Let $p$ be a prime dividing $h(L)$.
Let $r_p$ be the rank of the $p$-class group of $L$.
Then  $p\times \prod_{i=1}^{r_p}(p^i-1)$ and $N$ are not coprime.}
\item
We give also in this paper a connection with Geometry of Numbers point of view. With  an explicit  geometric  upper bound $H_F$ of the class number $h(F)$ for any  field $F$, which is  given in this paper: 

{\it Let $L/\Q$ be a finite Galois  extension with $[L:\Q]=N$. Let $h(L)$ be the class number of $L$. Suppose that $h(L)>1$. Suppose that $N$ has odd prime divisors. Let $q$ be an odd prime divisor of $N$.
Then  there exists a cyclic extension $L/F$ with  $q=[L:F]$. 
Suppose that  $p>H_F$ is  a prime dividing  $h(L)$. 
Let $r_p$ be the rank of the $p$-class group of $L$.
Then  $p\times (p^{r_p}-1)\equiv 0\modu q$.}
\item
The proofs are {\bf  elementary}.
We give  several verifications of results obtained for cyclic and abelian  extensions from the tables in  Washington\cite{was}, Schoof \cite{sch}, Masley \cite{mas}, Girstmair \cite{gir}, Jeannin \cite{jea} and the tables of cubic totally real number fields of the $ftp$ server
 {\it megrez.math.u-bordeaux.fr}.
\endabstract
\end{itemize}
%%% ====================================================================
%
%RRRRRRRR 10
%%% ====================================================================
%%% ====================================================================
%
%RRRRRRRR 10
%%% ====================================================================
\clearpage
\section{ On prime factors of  class number of finite algebraic extensions K/\Q}
The so called {\it rank theorem}, see Masley \cite{mas}, corollary 2.15 p. 305 is:

Suppose $M/P$ is a cyclic extension of degree $m$. Let $p$ be a prime which does not divide $h(E)$, class number of $E$  for 
any field $E$ with $P\subset E\subset M,\quad E\not=M$, and which does not divide $m$. If $p\ |\ h(M)$ class number of $M$ then the rank $r_p$ of $p$-class group of $M$ is a multiple of $f$, the order of $p\modu m$. 

In this subsection we use this theorem to get some congruences on prime factors of class number of some Galois extensions $K/\Q$.
 
\subsection{ Some definitions}
\begin{itemize}
\item
Let $K/\Q$ be a finite algebraic extension.
Let $L/K$ be a finite  Galois  extension. We say that $L/K$ is a {\it cyclic tower extension} 
if there is a tower of fields $K=L_0\subset L_1\subset \dots\subset L_n\subset L=L_{n+1}$ where 
$L_{j+1}/L_j,\quad j=0,\dots, n,$ is a cyclic extension. Observe that, without loss of generality, we can assume that $[L_{j+1}:L_j]$ is a prime.
\item
We say that $L/K$ is a Galois solvable extension 
if $L/K$ is a Galois extension with a solvable Galois group $Gal(L/K)$. 
\item
Let $L/K$ be a finite extension. We say that $L/K$ is a solvable extension if 
the smallest Galois extension $M$ of $K$ containing $L$ is Galois solvable, see for instance Lang, \cite{lan} paragraph 7 p. 216. Classically, the Galois solvable extension $M/K$ is a cyclic tower extension.
\item
Let $L/\Q$ be a finite  algebraic extension. Let  $N=[L:\Q]$,  
where $N=2^{\alpha_0}\times N_1$ with $N_1$  odd and where we assume that $N_1>1$. 
\item
Let  $h(L)$ be the class number of the field $L$. In
this paper, we  are studying, for fields $L$ with $h(L)>1$,  some congruences on the primes dividing $h(L)$.
\item
The class group $C$ of $L$ is an abelian group. Let $p$ be a prime dividing $h(L)$. The $p$-subgroup $C_p$ of $C$ is an abelian group,  direct sum of $r_p$ cyclic group whose order are power of $p$. 
This number $r_p$ is the $p$-rank of $p$-subgroup of $C$. Observe that if $p^{e_p}\| h(L)$ then $r_p\leq e_p$ and that  
$r_p$ is generally {\it small} and more often $r_p=1$.
\end{itemize}
%%% ====================================================================
%
%RRRRRRRR 10
%%% ====================================================================
%%% ====================================================================
%
%RRRRRRRR 10
%%% ====================================================================
\subsection{Some results}
%%% ====================================================================
%
%RRRRRRRR 10
%%% ====================================================================
%%% ====================================================================
%
%RRRRRRRR 10
%%% ====================================================================
\begin{thm}\label{t211291} 
Let $L/\Q$ be a finite  algebraic  extension   with $[L:\Q]=2^{\alpha_0}\times N_1$, where $N_1>1$ is odd.
Suppose that there exists a  field $K\subset L$ with $[K:\Q]=2^{\alpha_0}$ and with
$L/K$  Galois solvable extension.
Let $h(L)$ be the class number of $L$. Suppose that $h(L)>1$.
Let $p$ be a prime dividing $h(L)$.
Let $r_p$ be the rank of the $p$-class group of $L$.
If $p\times \prod_{i=1}^{r_p}(p^i-1)$ and $N_1$ are coprime,  then $p$ divides the class number $h(K)$ of $K$.
\begin{proof}$ $
\begin{itemize}
\item
$L/K$ being Galois solvable, we can always suppose that  
there exists a tower of  fields
\begin{displaymath}
K\subset L_1\subset L_2\subset \dots\subset L_j\subset L_{j+1}\subset\dots\subset L_t\subset L,
\end{displaymath}
such that 
$[L_{j+1}:L_j]=q_j$ where $L_{j+1}/L_j,\quad j=1,\dots,t$, is a cyclic extension with $q_j$  prime.
\item
The extension $L/L_t$ is cyclic with $[L:L_t]=q_t$ and $q_t|N_1$.
\item
Suppose at first that $p\not|h(L_t)$:
then hypotheses of   rank theorem are verified, see Masley \cite{mas} corollary 2.15 p. 305 : 
\begin{itemize}
\item
From hypothesis,  $p$ does not divide $N_1$, so $p$ does not divide $q_t=[L:L_t]$.
\item
The extension $L/L_t$ is cyclic with 
$ h(L)\equiv 0\modu p,\quad  h(L_t)\not\equiv 0\modu p$.
\item
There is no field $E$, different of $L$ and of $L_t$ with $L_t\subset E\subset L$.
\end{itemize}
From rank theorem, if $f$ is the order of $p \modu q_t$ then $f|r_p$ and so $p^{r_p}-1\equiv 0\modu q_t$,
which contradicts hypothesis. 
\item
Therefore   $p$ divides $h(L_t)$ : the $p$-rank $r_p^\prime$  of class group $h(L_t)$ of $L_t$ verifies $r_p^\prime\leq r_p$ because $h(L_t)\ |\ q_t\times h(L)$,
 see for instance Masley, \cite{mas}, {\it pushing up} corollary 2.2 p. 301 
and so $p$ and $q_t$  coprime   implies that the $p$-class group $C_p(L_t)$ of $L_t$ and $C_p$ of $L$ verify 
$C_p(L_t)\subset C_p(L)$ and so $r_p^\prime\leq r_p$. Then  let us consider 
$N^\prime=\frac{N}{q_t}=[L_t:K]$. 
We can pursue the same algorithm with same Masley corollaries 2.15 and 2.2 applied to  extension $L_t/K$  in place of  $L/K$  
and $N^\prime$ in place of $N$ and with $r_p^\prime\leq r_p$ in place of $r_p$,  up to find a prime divisor 
$q_j,\quad 1\leq j\leq t,$  of $N_1$ dividing $\prod_{i=1}^{r_p} (p^i-1)$,  or to get the  subfield   $K$ of $L$ with $[K:\Q]=2^{\alpha_0}$ and $p |h(K)$ , which achieves the proof.
\end{itemize}
\end{proof}
\end{thm}
%%% ====================================================================
%
%RRRRRRRR 10
%%% ====================================================================
In following corollary, we generalize the result obtained for  $L/K$  Galois solvable extension  to the case where $L/K$ is a solvable extension.
\begin{cor}\label{c302151} 
Let $L/\Q$ be a finite  algebraic  extension   with $[L:\Q]=2^{\alpha_0}\times N_1$, where $N_1>1$ is odd.
Suppose that there exists a  field $K\subset L$ with $[K:\Q]=2^{\alpha_0}$ and with
$L/K$   solvable extension.
Let $M/K$ be the Galois solvable extension containing $L/K$.
Let $h(L)$ be the class number of $L$. Suppose that $h(L)>1$.
Let $p$ be a prime dividing $h(L)$ such that  $p>[L:K]$.
Let $r_p(M)$ be the rank of the $p$-class group of $M$.
If $\prod_{i=1}^{r_p(M)}(p^i-1)$ and $N_1$ are coprime,  then $p$ divides the class number $h(K)$ of $K$.
\begin{proof}$ $
From {\it pushing up} theorem, see Masley \cite{mas} corollary 2.2 p. 301, $h(L)\ |\ [M:L]\times h(M)$.
Let  $q_M$ be the greatest prime divisor of $[M:K]$. 
The solvable extension $M/K$ is the smallest Galois extension of $K$ containing $L$. Therefore 
$q_M\  |\  [L:K]!$, see for instance Morandi, \cite{mor}, cor. 3.8  p. 29. From hypothesis, this implies that $p\not |\  [M:L]$. Therefore 
$p\ |\ h(M)$. Then apply theorem \ref{t211291} p. \pageref{t211291} to Galois solvable extension $M/K$.
\end{proof}
\end{cor}
%%% ====================================================================
%
%RRRRRRRR 10
%%% ====================================================================
In following corollary  we give a particularly straightforward formulation 
of theorem \ref{t211291} p. \pageref{t211291} when $[L:\Q]=N$ is odd. 
\begin{cor}\label{c302112} 
Let $L/\Q$ be a   Galois solvable   extension   with $[L:\Q]= N$, where $N>1$ is odd.
Let $h(L)$ be the class number of $L$. Suppose that $h(L)>1$.
Let $p$ be a prime dividing $h(L)$.
Let $r_p$ be the rank of the $p$-class group of $L$.
Then  $p\times \prod_{i=1}^{r_p}(p^i-1)$ and $N$ are not coprime.
\begin{proof}$ $
Apply theorem \ref{t211291} p. \pageref{t211291}.
\end{proof}
\end{cor}
%%% ====================================================================
%
%RRRRRRRR 10
%%% ====================================================================
\paragraph{Remark:} Observe that, with previous notations,  if $L/K$ is abelian then $L/K$ is Galois solvable and so 
the previous results can be applied to all abelian extensions $L/K$ with $[K:\Q]=2^{\alpha_0}$ and 
$[L:K]$ odd.
%%% ====================================================================
%
%RRRRRRRR 10
%%% ====================================================================
\section{Geometry of Numbers point of view}
\begin{itemize}
\item
The next result connects rank theorem to Geometry of Numbers point of view.
Let $L/\Q$ be a finite  algebraic extension.  
Suppose  that there exists a subfield $F$ of $L$ such that $L/F$ is a cyclic extension, with $[L:\Q]=N,\quad [L:F]=q$ where $q$ is an odd prime.
\item
Let   $D_F$ be the discriminant of $F$. Let $m=\frac{N}{q}$.
Then  $h(F)\leq \frac{2^{m-1}}{(m-1)!}\times \sqrt{|D_F|}\times (log(|D_F|))^{m-1}$, 
see Bordell\`es \cite{bor} theorem 5.3 p. 4.
Let us note 
\begin{equation}\label{e212141}
H_F=\frac{2^{m-1}}{(m-1)!}\times \sqrt{|D_F|}\times (log(|D_F|))^{m-1}. 
\end{equation}
If a prime $p$ verifies
$p>H_F$ then $p\not| h(F)$. 
\end{itemize}
%%% ====================================================================
%
%RRRRRRRR 10
%%% ====================================================================
\begin{thm}\label{t212141}
Let $L/\Q$ be a finite algebraic extension. Let $h(L)$ be the class number of $L$. Suppose that $h(L)>1$.
Suppose that there exists a cyclic extension $L/F$, where $q=[L:F]$ is an odd prime.
Let $H_F$ be a geometric upper bound of class number of $F$ given by relation (\ref{e212141}).
Suppose that  $p>H_F$ is  a prime dividing  $h(L)$. 
Let $r_p$ be the rank of the $p$-class group of $L$.
Then  $p\times (p^{r_p}-1)\equiv 0\modu q$.
\begin{proof}
We apply rank theorem, see Masley\cite{mas} corollary 2.15 p 305, and upper bound
$H_F$ of class number $h(F)$ of field $F$ given in relation (\ref{e212141}). 
\end{proof}
\end{thm}
%%% ====================================================================
%
%RRRRRRRR 10
%%% ====================================================================
\begin{cor}\label{c212151}
Let $L/\Q$ be a finite Galois  extension with $[L:\Q]=N$. Let $h(L)$ be the class number of $L$. Suppose that $h(L)>1$. Suppose that $N$ has odd prime divisors. Let $q$ be an odd prime divisor of $N$.
Then  there exists a cyclic extension $L/F$ with  $q=[L:F]$. 
Let $H_F$ be a geometric upper bound of class number of $F$ given by relation (\ref{e212141}).
Suppose that  $p>H_F$ is  a prime dividing  $h(L)$. 
Let $r_p$ be the rank of the $p$-class group of $L$.
Then  $p\times (p^{r_p}-1)\equiv 0\modu q$.
\begin{proof}
Immediate consequence of theorem \ref{t212141} and of Galois theory. 
\end{proof}
\end{cor}
\paragraph{Remark:} Observe that all the results above  give strong information on prime factor of class numbers of Galois solvable extensions of $L/K$, because practically the $p$-rank $r_p$ is generally {\it small}: in all numerical examples following $r_p\leq 3$ and it verifies  more often $r_p=1$. 
%%% ====================================================================
%
%RRRRRRRR 10
%%% ====================================================================
%%% ====================================================================
%
%RRRRRRRR 10
%%% ====================================================================
\section{Numerical examples}\label{ss22043}
The examples found to  check theses results are taken from: 
\begin{itemize} 
\item
the table of relative class numbers of cyclotomic number fields in Washington, \cite{was} p 412,  with some elementary MAPLE computations,
\item 
the table of relative class number of cyclotomic number fields in Schoof, \cite{sch}
\item
the table of maximal real subfields $\Q(\zeta_l+\zeta_l^{-1})$ of $\Q(\zeta_l)$ for $l$ prime in Washington, \cite{was} p 420.
\item 
the tables of relative class number of imaginary cyclic fields of Girstmair of degree 4,6,8,10. \cite{gir}.
\item
the tables of quintic number fields computed by Jeannin, \cite{jea}
\item 
the tables of cubic totally real cyclic number fields of the Bordeaux University in the Server 

{\it megrez.math.u-bordeaux.fr}.
\end{itemize}

All the results examined in these tables are in accordance with our theorems.
%%% ====================================================================
%
%RRRRRRRR 10
%%% ====================================================================
\subsection{Cyclotomic number fields $\Q(\zeta_u)$}
 
Let $u\in\N,\quad u>2$. Let $\zeta_u$ be a primitive $u^{th}$ root of unity.
Here we have $N=\phi(u)$, where $\phi$ is the Euler indicator.
The cyclotomic number fields $\Q(\zeta_u)$ of the examples are taken  with $2\| N=\phi(u)$,
except the example $\zeta_u$ with $u=572,\quad N=\phi(u)=2^4.3.5$.
 
\begin{itemize}
\item 

$\Q(\zeta_u), \quad u=59 ,\quad N=\phi(u)=58=2.29 :$

$h^- = 3 . (2.29+1). (2^3 . 29 +1),\quad  h^+=1.$

$3 =h(\Q(\sqrt{-59}).$
\item

$\Q(\zeta_u), \quad u=71, \quad N=\phi(u)=70=2.5.7$ :

$h^- = 7^2.(2^3.5.7.283+1), \quad h^+=1. $

\item
$\Q(\zeta_u), \quad u=79,\quad  N=\phi(u)=78=2.3.13$ :

$ h^-= 5.(2^2.13+1)(2.3^2.5.13.17.19+1),\quad h^+=1.$

$5=h(\Q(\sqrt{-79})$.

\item
$\Q(\zeta_u), \quad u=83,\quad N=\phi(u)=82=2.41$ :

$ h^-=3.(2^2.41.1703693+1),\quad h^+=1$.

$3=h(\Q(\sqrt{-83}).$
\item
$\Q(\zeta_u),\quad u=103,\quad  N=\phi(u)=102=2.3.17$:

$h^-= 5.(2.3.17+1)(2^2.3.5.17+1)(2.3^2.5.17.11273+1)$,

$h^+=1$.

\item
$\Q(\zeta_u), \quad u=131,\quad N=\phi(u)=130=2.5.13$ :

$ h^-=3^3.5^2.(2^2.13+1)(2.5.13+1)(2^2.5^2.13+1)\times$

$(2^2.5^2.13.29.151.821+1), \quad h^+=1$.

Observe that for $p=3$, we have  $p^3\|h$ and the group $C_p$ is not cyclic as it is seen in Schoof \cite {sch} table 4.2 p 1239 , where $r_p=3$. The rank theorem in that case shows that $f=r_p=3$ and so that $3^f=27\equiv 1 \modu 13$.

\item
$\Q(\zeta_u), \quad u=151,\quad N=\phi(u)=150=2.3.5^2$ :

$ h^-=(2.3+1)(2.5+1)^2.(2^2.5.7+1)(2.3.5^2.173+1)\times$

$(2^2.3.5^4.7.23+1)(2^2.3.5^2.7.13.73.1571+1),\quad h^+=1.$

For $p=2.5+1=11$, from Schoof, see table 4.2 p 1239,  the group is not cyclic, $r_p=2$. From rank theorem, we can only assert that $f|r_p=2$, so that $f=1$ or $f=2$, so only that $p^2\equiv 1 \modu g=5$. We see in this numerical application that in that case  $f=1$ and $f<r_p$ strictly.

\item $\Q(\zeta_u), \quad u=191,\quad N=\phi(u)=190=2.5.19$:

$ h^-= (2.5+1).13.(2.19^2.71+1)(2.3.5.19.277.3881+1),$

$h^+= (2.5+1)$.

$13=h(\Q(\sqrt{-191})$.

Here, we note that the prime $p=11$ corresponds to $11$-class group of   
$\Q(\zeta_{191}+\zeta_{191}^{-1})$.

\item
$\Q(\zeta_u),\quad u=572=2^2.11.13,\quad N=\phi(u)=2^4.3.5,$

$h^-=3.5^2.(2.3+1)(2.3^2+1)^2(2.3.5+1)(2^3.5+1)(2^2.3.5+1)^2$

$(2^2.3.5.7+1)(2^2.3.5.11+1)(2.3^2.5.307+1)(2.3^2.5.11.73+1)$

$(2^2.3^2.5^2.31+1)(2^2.3^5.5.7.53.263+1)$

We don't know $h^+$.

\end{itemize}
{\bf Remarks:} 
\begin{itemize}
\item
Let $[K:\Q]=2^{\alpha_0}\times n_1^{\alpha_1}\times\dots\times n_k^{\alpha_k}$. We observe in  Washington  \cite{was}, tables  of relative class numbers  p. 412 and of real class numbers p. 421  that, when $p$ is large, $p-1$ is divisible by several  or all primes in the set 
$\{n_i\ | \  i=1,\dots,k\}$. This observation complies with Geometry of Number corollary 
\ref{c212151} p. \pageref{c212151}.
\item
Observe that frequently , we get $f=r_p=1$ and so $p\times (p-1)\equiv 1\modu n$. In our examples we get 
only one example $\Q(\zeta_{131})$ with $p=3,\quad f=r_p=3$. 
\end{itemize}
%%% ====================================================================
%
%RRRRRRRR 10
%%% ====================================================================
\subsection{Real class number $\Q(\zeta_l+\zeta^{-1}_l)$}
The examples are obtained from the table of real class number in Washington, \cite{was} p 420.
Here, $l$ is a prime, the class number  $h_\delta$ is the conjectured value of the class number $h^+$ of 
$\Q(\zeta_l+\zeta_l^{-1})/\Q$ with a minor incertitude on an extra factor. But the factor $h_\delta$ must verify our theorems. 
We extract some examples of the table with $2\| l-1$. Here we have $N=\frac{l-1}{2}$ and thus 
$N\not \equiv 0 \modu 2$. 
\begin{itemize}
\item
$l=191,\quad l-1=2.5.19,\quad h_\delta=(2.5+1)$

\item
$l=1063,\quad l-1=2.3^2.59,\quad h_\delta=(2.3+1).$

\item
$l=1231,\quad l-1=2.3.5.41,\quad h_\delta=(2.3.5.7+1).$

\item
$l=1459,\quad l-1=2.3^6,\quad h_\delta=(2.3.41+1).$

\item
$l=1567,\quad l-1=2.3^3.29,\quad h_\delta=(2.3+1).$

\item
$l=2659,\quad l-1=2.3.443,\quad h_\delta=(2.3^2+1).$

\item
$l=3547,\quad l-1=2.3^2.197,\quad h_\delta=((2.3^2+1)(2.3^2.7^2+1).$

\item
$l=8017,\quad l-1=2^4.3.167$,

$h_\delta=(2.3^2+1)(2.3^2.7^2+1)(2^2.3^3+1).$

\item
$l=8563,\quad l-1=2.3.1427,\quad h_\delta=(2.3+1)^2$. We have $r_p=1$ or $r_p=2$ so $f=1$ or $f=2$. We can conclude from theorem \ref{t211291} p. \pageref{t211291}, that $p^2\equiv 1\modu 3$. 

\item
$l=9907,\quad l-1=2.3.13.127,\quad h_\delta=(2.3.5+1).$
\end{itemize}

%%% ====================================================================
%
%RRRRRRRR 10
%%% ====================================================================

\subsection {Cubic fields $K/\Q$ cyclic and totally real.}
Here, we have $N=3$. Note that in that case discriminants are square in $\N$. 
\begin{itemize}
\item Discriminant $D_1=3969=(3^2.7)^2$, \quad$h=3$
\item Discriminant $D_2=3969=(3^2.7)^2$, \quad $h=7=2.3+1$
\item Discriminant $D_1=8281=(7.13)^2$, \quad $h=3$
\end{itemize}
%\end{itemize}
%%% ====================================================================
%
%RRRRRRRR 10
%%% ====================================================================
\subsection{Totally real cyclic fields of prime conductor $<100$}
We have found few numeric results in the literature. We refer to Masley, \cite{mas}, Table 3 p 316.
In these results, $K/\Q$ is a real cyclic field with $[K:\Q]=N$, with conductor $f$, with root of discrimant $Rd$ and class number $h$.
\begin{itemize}
\item 
$f=63,\quad N=3,\quad Rd=15.84,\quad \quad h=3$
\item 
$f=63,\quad N=3,\quad Rd=15.84,\quad \quad h=3$
\item 
$f=63,\quad N=6,\quad Rd=26.30,\quad \quad h=3$
\item 
$f=63,\quad N=6,\quad Rd=26.30,\quad \quad h=3$
\item 
$f=91,\quad N=3,\quad Rd=20.24,\quad \quad h=3$
\end{itemize}
%%% ====================================================================
%
%RRRRRRRR 10
%%% ====================================================================
\subsection{Lehmer quintic cyclic field}
The prime divisors of the $82$ cyclic number fields of the table in Jeannnin \cite{jea}, with conductor $f<3000000$,  are, at a glance, of the form $p=2$ or $p=5$ or $p\equiv 1 \modu 10$, which clearly verifies corollary \ref{c302112} p.\pageref{c302112}. 
%%% ====================================================================
%
%RRRRRRRR 10
%%% ====================================================================
\subsection{Decimic imaginary cyclic number fields with conductor between  $9000$ and $9500$}
This example is obtained from the tables of Girstmair, \cite{gir}. $f$ is a prime conductor, $h$ is the factorization of the class number $K/\Q$.
\begin{itemize}
\item 
$f=9011$, $h=3.(2.5+1).(2.3.5.52201+1)$.

$3$ divides the class number of $\Q(\sqrt{-9011})$.
\item 
$f=9081$, $h=3.7.(2^3.5+1)$.

$3,7$ divide the class number of $\Q(\sqrt{-9081})$.
\item 
$f=9151$, $h=67.(2^2.5.1187+1)$.

$67$ divides the class number of $\Q(\sqrt{-9151})$.
\item 
$f=9311$, $h=97.(2.5.5689+1)$.

$97$ divides the class number of $\Q(\sqrt{-9311})$.
\item 
$f=9371$, $h=7^2.(2.3.5^2+1).(2.3^3.5+1)$.

$7$ divides the class number of $\Q(\sqrt{-9371})$.
\item 
$f=9391$, $h=5^2.(2^4.3.5.7.71+1)$.
\end{itemize} 
%%% ====================================================================
%
%RRRRRRRR 10
%%% ====================================================================
\subsection{Not abelian Galois solvable extensions $L/K$ }
Here,  we give an example of Galois solvable extension not abelian.
Let $q$ be a prime with $q=2g+1$, where $g$ is an odd prime.
Let $\zeta_q$ be a  root of the equation $X^{q-1}+\dots+X+1=0$ and let $\omega$ root of the equation 
$X^q-2=0$.
The extension $L=\Q(\zeta_q,\omega)/\Q=\Q(\zeta_q)(\omega)$ verifies 
$\Q\subset\Q(\sqrt{-q})\subset\Q(\zeta_q)\subset L$, where $\Q(\zeta_q)/\Q(\sqrt{-q})$ is cyclic and where 
$\Q(\zeta_q,\omega)/\Q(\zeta_q)$ is cyclic. Let  $p$ be  a prime dividing 
the class number $h(L)$ with $ p\not=g,\quad p\not= q$. 
Apply theorem \ref{t211291} p. \pageref{t211291} to get:

Let $r_p$ be the rank of the $p$-class group of $L$.
Then else $p\ |\ h(\Q(\sqrt{-q}))$ else $Gcd(\prod_{i=1}^{r_p}(p^i-1), q\times g)>1$.
%%% ====================================================================
%\maketitle
%
%RRRRRRRR 10
%%% ====================================================================
%%% ====================================================================
%
%RRRRRRRR 10
%%% ====================================================================
%
%RRRRRRRR 560
%%% ====================================================================

***************

Roland Qu\^eme

13 avenue du ch\^ateau d'eau

31490 Brax

France

2003 april 20

tel 0561067020

e-mail : roland.queme@free.fr

home page: http://roland.queme.free.fr/

\end{document}